\documentclass{amsart}

\usepackage{mathtools}
\usepackage{amsmath}
\usepackage{amssymb}

\newtheorem{theorem}{Theorem}[section]
\newtheorem{lemma}[theorem]{Lemma}
\newtheorem{proposition}[theorem]{Proposition}
\newtheorem{corollary}[theorem]{Corollary}

\newtheorem*{theorem*}{Theorem}

\theoremstyle{definition}
\newtheorem{definition}[theorem]{Definition}

\newtheorem*{definition*}{Definition}

\theoremstyle{remark}
\newtheorem{remark}[theorem]{Remark}

\numberwithin{equation}{section}

\newcommand{\cont}{{\mathrm C }}

\newcommand{\catb}[1]{{\mathrm B}(\ca{#1}) }
\newcommand{\catl}[1]{{\mathrm L}^\infty (\ca{#1},\meagre{#1}) }
\newcommand{\catn}[1]{{\mathrm N}(\ca{#1},\meagre{#1}) }

\newcommand{\contK}{\cont(K)}

\newcommand{\contKd}{\contK^\delta}

\newcommand{\cfam}[1]{{\mathrm M} (\ca{#1},\meagre{#1}) }

\newcommand{\borm}{\mathrm M}
\newcommand{\bormK}{\borm (K)}
\newcommand{\famr}[1]{\mathrm M (\roa{#1})}

\newcommand{\zerof}{{\mathbf 0}}
\newcommand{\onef}{{\mathbf 1}}

\newcommand{\N}{\mathbb{N}}

\newcommand{\R}{\mathbb{R}}

\newcommand{\roa}[1]{\mathcal{R}_{#1}}
\newcommand{\ca}[1]{\Sigma^\mathcal{C}_{#1}}
\newcommand{\bor}[1]{\Sigma^\mathcal{B}_{#1}}
\newcommand{\meagre}[1]{\mathcal{N}_{#1}}

\usepackage[UKenglish]{babel}

\usepackage{currfile}
\usepackage[UKenglish]{datetime}
\usepackage{color}

\begin{document}

\title{Category measures, the dual of $\contKd$ and hyper-Stonean spaces}

\author{Jan Harm van der Walt}
\address{Department of Mathematics and Applied Mathematics, University of Pretoria, Cor\-ner of Lynnwood Road and Roper Street,
Hatfield 0083, Pretoria, South Africa}



\begin{abstract}
For a compact Hausdorff space $K$, we give descriptions of the dual of $\contKd$, the Dedekind completion of the Banach lattice $\contK$ of continuous, real-valued functions on $K$.  We characterize those functionals which are $\sigma$-order continuous and order continuous, respectively, in terms of Oxtoby's category measures.  This leads to a purely topological characterization of hyper-Stonean spaces.
\end{abstract}

\maketitle

\section{Introduction}\label{Section:  Introduction}

It is well known that the Banach lattice $\contK$ of continuous, real-valued functions on a compact Hausdorff space $K$ fails, in general, to be Dedekind complete.  In fact, $\contK$ is Dedekind complete if and only if $K$ is Stonean, i.e. the closure of every open set is open, see for instance \cite{Schaefer1974}.  A natural question is to describe the Dedekind completion of $\contK$; that is, the unique (up to linear lattice isometry) Dedekind complete Banach lattice containing $\contK$ as a majorizing, order dense sublattice.  Several answers to this question have been given, see for instance \cite{Anguelov2004,deJongevanRooijRieszSpaces1977,Dilworth1950,Horn1953,Maharam1976,Maxey,vanderwaltPNA2019} among others.  A further question arises: what is the dual of $\contKd$?  As far as we are aware, no direct answer to this question has been stated explicitly in the literature.

We give two answers to this question; one in terms of Borel measures on the Gleason cover of $K$ \cite{Gleason1958}, and another in terms of measures on the category algebra of $K$ \cite{Oxtoby1960,Oxtoby1980}.  We use the latter description to obtain natural characterizations of strictly positive, $\sigma$-order continuous, and order continuous functionals on $\contKd$.  In particular, we show that $\contKd$ admits a strictly positive, $\sigma$-order continuous functional if and only if $K$ admits a category measure in the sense of Oxtoby \cite{Oxtoby1960}.  Specialising to the case of a Stonean space $K$, we obtain a purely topological characterization of hyper-Stonean spaces.  Although the class of hyper-Stonean spaces was introduced by Dixmier \cite{Dixmier1951}  in 1951, no such characterization has been given to date \cite[page 197]{DalesDashiellLauStrass2016}.

The paper is organised as follows.  In Section \ref{Section:  Preliminaries} we introduce notation and recall definitions and results from the literature to be used throughout the rest of the paper.  Section \ref{Section:  Characterisations of C(K)^d} contains two representations of $\contKd$ which are amendable to the problem of characterizing the dual of this space; these characterizations of $(\contKd)^\ast$ are given in Section \ref{Section:  The dual of C(K)^d}.  Strictly positive, $\sigma$-order continuous and order continuous functionals on $\contKd$ are discussed in Section \ref{Section:  Order continuous elements of (C(K)^d)^*} which ends with a characterization of hyper-Stonean spaces.

\section{Preliminaries}\label{Section:  Preliminaries}

\subsection{Vector and Banach lattices}

Let $E$ be a real vector lattice. We assume that $E$ is Archimedean; that is, $\inf\{\frac{1}{n} f ~:~ n\in\N\}=0$ for every $f\geq 0$.  Denote by $E_+$ the positive cone in $E$, $E_+=\{f\in E ~:~ 0\leq f\}$. For $f\in E$ the positive part, negative part and modulus of $f$ is given by
\[
f^+ = f\wedge 0,~ f^- = (-f)\wedge 0 \text{ and } |f|=f^+ + f^-,
\]
respectively. For $D\subseteq E$ and $f\in E$ we write $D\downarrow f$ if $D$ is downward directed and $\inf D=f$.

$E$ is \emph{Dedekind complete} if every subset of $E$ which is bounded above (below) has a least upper bound (greatest lower bound).  We call $E$ \emph{order separable} if every subset $A$ of $E$ contains a countable set with the same set of upper bounds as $A$.

Let $F$ be a vector lattice subspace of $E$.  $F$ is \emph{order dense} in $E$ if for every $0<f\in E$ there exists $g\in F$ so that $0<g\leq f$.  If for every $f\in E_+$ there exists $g\in F$ so that $f\leq g$ then $F$ is \emph{majorizing} in $E$.  An \emph{ideal} in $E$ is a vector lattice subspace $F$ of $E$ so that if $f\in \in F$ and $|g|\leq |f|$ then $g\in F$. A $\sigma$-ideal in $E$ is an ideal $F$ in $E$ so that $\sup A\in F$ whenever $A$ is a countable subset of $F$ and $\sup A$ exists in $E$.  The \emph{principal ideal} generated by $f\in E$ is $I_f = \{g\in E ~:~ |g|\leq \alpha |f| \textit{ for some }\alpha\in\R_+\}$.

An element $e$ of $E_+$ is called an \emph{order unit} if for every $f\in E_+$ there exits a real number $\alpha >0$ so that $f\leq \alpha e$.  A \emph{weak order unit} is any $e\in E_+$ so that if $f \wedge e = 0$ then $f=0$. Observe that $e\in E_+$ is an order unit if and only if $I_e=E$; $e$ is a weak order unit if and only if $I_e$ is order dense in $E$.

A linear functional $\varphi$ on $E$ is \emph{order bounded} if it maps order bounded subsets of $E$ to bounded subsets of $\R$, and \emph{positive} if $\varphi(f)\geq 0$ whenever $f\in E_+$.  Every positive functional on $E$ is order bounded.  Conversely, every order bounded functional is the difference of two positive functionals. The space $E^\sim$ of all order bounded functionals on $E$ is a Dedekind complete vector lattice with respect to the ordering
\[
\varphi\geq \psi \text{ if and only if } \varphi-\psi \text{ is positive}.
\]
Therefore every $\varphi\in E^\sim$ can be expressed as the difference of positive functionals, $\varphi = \varphi^+ - \varphi^-$.

If $E$ is Banach lattice then a functional $\varphi$ on $E$ is order bounded if and only if it is norm bounded.  Therefore the norm dual $E^\ast$ and order dual $E^\sim$ of $E$ coincide and is itself a Banach lattice.

A functional $\varphi\in E^\sim$ is \emph{$\sigma$-order continuous} if $\inf\{|\varphi(f_n)| ~:~ n\in\N\}=0$ whenever $f_n\downarrow 0$ in $E$. If $\inf\{|\varphi(f)| ~:~ f\in D\}=0$ for every $D\downarrow 0$ in $E$, we say that $\varphi$ is \emph{order continuous}.  A functional $\varphi$ is ($\sigma$-)order continuous if and only if $|\varphi|$ is ($\sigma$-)order continuous, if and only if both $\varphi^+$ and $\varphi^-$ are ($\sigma$-)order continuous.

\begin{lemma}\label{Lemma:  Restriction and extension of order continuous functionals}
Let $E$ be an Archimedean vector lattice and $F$ an order dense, majorizing vector lattice subspace of $E$. \begin{enumerate}
    \item[(i)] If $\varphi\in E^\sim$ then $\varphi$ is order continuous on $E$ if and only if the restriction of $\varphi$ to $F$ is order continuous.
    \item[(ii)] If $\varphi\in E^\sim$ is $\sigma$-order then the restriction of $\varphi$ to $F$ is $\sigma$-order continuous.
    \item[(iii)] If $\varphi\in F^\sim$ is order continuous, then $\varphi$ has a unique order continuous extension to $E$.
\end{enumerate}
\end{lemma}

\begin{proof}
We prove the results for positive functionals.  The case of a general $\varphi\in E^\sim$ follows from the decomposition $\varphi=\varphi^+-\varphi^-$.

\noindent {\it Proof of (i).}~ Let $D\subseteq F_+$.  Because $F$ is order dense in $E$, $D\downarrow 0$ in $E$ if and only if $D\downarrow 0$ in $F$. Therefore, if $\varphi$ is order continuous on $E$ then its restriction to $F$ is order continuous.

Conversely, assume that the restriction of $\varphi$ to $F$ is order continuous.  Let $D\downarrow 0$ in $E$.  For each $g\in D$ let $D_g = \{f\in F ~:~ g\leq f\}$, and let $D'=\bigcup\{D_g ~:~g\in D\}$.  The fact that $F$ is majorizing in $E$ guarantees that $D_g\neq \emptyset$ for each $g\in D$.  We observe that $D\downarrow 0$ in $E$ implies that $D'$ is downward directed. In particular, if $f,h\in D'$ then $f\wedge h\in D'$. Because $F$ is an order dense sublattice of $E$, $g=\inf D_g$ for every $g\in D$.  Therefore $D'\downarrow 0$ in $E$, hence also in $F$.  It follows from the positivity and order continuity of $\varphi$ on $F$ that $0\leq \inf \varphi[D] \leq \inf \varphi[D']=0$. Therefore $\inf \varphi[D]=0$.

\noindent {\it Proof of (ii).}~ The proof is similar to the proof of (i).

\noindent {\it Proof of (iii).}~ This is a special case of \cite[Theorem 4.12]{AliprantisBurkinshaw2006}.
\end{proof}

\subsection{$\contK$ and its dual}

Let $K$ be a compact Hausdorff space.  $\contK$ is the real Banach lattice of continuous, real-valued functions on $K$; the algebraic operations and order relation are defined pointwise as always, and the norm on $\contK$ is $\|f\|=\max\{f(x)| ~:~ x\in K\}$.  The symbols $\zerof$ and $\onef$ denote the constant functions with values $0$ and $1$, respectively.  For $A\subseteq K$, the characteristic function of $A$ is denoted $\onef_A$.

The space $K$ satisfies the countable chain condition (ccc) if every pairwise disjoint collection of nonempty open sets in $K$ is countable.  $\contK$ is order separable if and only if $K$ satisfies ccc \cite[Theorem 10.3]{dePagterHuijsmans1980II}.

$\bor{K}$ denotes the $\sigma$-algebra of Borel sets in $K$ and $\bormK$ the space of regular signed Borel measures on $K$; that is, $\sigma$-additive, real valued regular measures on $\bor{K}$.  As is well known, $\bormK$ is a Banach lattice with respect to the total variation norm
\[
\|\varphi\| = |\varphi|(K),~ \varphi\in \bormK
\]
and pointwise ordering.  In fact, $\bormK$ is isometrically lattice isomorphic to $\contK^\ast$.  We will not make a notational distinction between a measure $\varphi$ and the associated functional.  A measure $\varphi\in \bormK$ is called \emph{normal} \cite{DalesDashiellLauStrass2016,Dixmier1951} if the associated functional on $\contK$ is order continuous. We denote the set of normal measures on by $\mathrm{N}(K)$.  Normal measures can be characterised as follows: $\varphi\in\bormK$ is normal if and only if $|\varphi|(N)=0$ for every meagre set $N\in \bor{K}$, see for instance \cite[Theorem 7.4.7]{DalesDashiellLauStrass2016}.

Recall that $K$ is \emph{Stonean} if the closure of every open set is open.  As mentioned in the introduction, $K$ is Stonean if and only if $\contK$ is Dedekind complete.  If $K$ is Stonean and, in addition, the union of the suppports of the normal measures on $K$ is dense in $K$, we call $K$ \emph{hyper-Stonean}, see \cite{Dixmier1951} and \cite[Section 4.7]{DalesDashiellLauStrass2016}.   We recall the following result which contains the main motivation for studying hyper-Stonean spaces, see for instance \cite{DalesDashiellLauStrass2016}.

\begin{theorem}
A space $K$ is hyper-Stonean if and only if any one of the following conditions hold. \begin{enumerate}
    \item[(i)] For every $\zerof<f\in\contK$ there exists a positive order continuous functional $\varphi$ on $\contK$ so that $\varphi(f)>0$.
    \item[(ii)] $\contK$ is isometrically isomorphic to the dual of a Banach space.
    \item[(iii)] $\contK$ is isometrically lattice isomorphic to the dual of a Banach lattice.
    \item[(iv)] The $C^\ast$-algebra of complex-valued continuous functions on $K$ is a von Neumann algebra.
\end{enumerate}
\end{theorem}

\subsection{Category measures}\label{Subsection:  Category measures}

Recall that a subset $A$ of $K$ has the property of Baire if there exists an open set $U$ and a meagre set $N$ in $K$ so that $A=U\Delta N$. Let
\begin{eqnarray}
\begin{array}{c}
\meagre{K}=\{N\subset K~:~N~{\rm is~meagre}\},\medskip\\
\ca{K} = \{U\Delta N~:~ U\text{ is open in }K,~N\in\mathcal{N}\}.\\
\end{array}\label{EQ:  Definitions SigmaC, N}
\end{eqnarray} 
The collection $\ca{K}$ is a $\sigma$-algebra \cite[Theorrem 4.3]{Oxtoby1980} and $\meagre{K}$ is a $\sigma$-ideal in $\ca{K}$.  In fact, $\meagre{K}$ is a $\sigma$-ideal in the powerset of $K$. The quotient algebra $\ca{K}/\meagre{K}$ is a complete Boolean algebra.  In fact, $\ca{K}/\meagre{K}$ is isomorphic to the algebra $\roa{K}$ if regular open subsets of $K$, see \cite{Maharam1976,Oxtoby1960}.

A \emph{category measure} on $K$ is positive, $\sigma$-additive measure $\mu:\ca{K}\to \R$ with the property that for any $A\in\ca{K}$, $\mu(A)=0$ if and only if $A\in\meagre{K}$, see for instance \cite{Oxtoby1960}.  There is a bijective correspondence between category measures on $K$ and strictly positive, $\sigma$-additive measures on $\roa{K}$. Oxtoby \cite{Oxtoby1960,Oxtoby1980} investigated the problem of characterizing those topological spaces which admit a category measure.  A solution to this problem was given by J. M. Ayerbe Toledano \cite{Ayerbe_Toledano1990}\footnote{We formulate the result for compact Hausdorff spaces, but it is valid for arbitrary Baire spaces.}.  In order to formulate this result we recall the following definitions, see \cite{Argyros1983,Kelley1959,Ayerbe_Toledano1990}.

\begin{definition}\label{Definition:  Intersection number}
Let $X$ be a nonempty set and $\mathcal{F}$ a collection of nonempty subsets of $X$.  For any $n\in\N$ and $\bar F = \langle F_1,\dots,F_n\rangle \in \mathcal{F}^n$ let
\[
i(\bar F) = \max \{ |J| ~:~ J\subseteq \{1,\ldots n\},~ \bigcap_{j\in J}F_j\neq \emptyset \}.
\]
The Kelley intersection number of $\mathcal{F}$ is defined as
\[
k(\mathcal{F}) = \inf \{\frac{i(\bar F)}{n} ~:~ n\in\N,~ \bar F \in \mathcal{F}^n \}.
\]
\end{definition}

\begin{definition}\label{Definition:  Property (***)}
Let $K$ be a compact Hausdorff space and $\mathcal{T}$ the collection of nonempty, open subsets of $K$.  We say that $K$ satisfies property (***) if there exists a partition $\{\mathcal{T}_n ~:~ n\in\N\}$ of $\mathcal{T}$ such that the following conditions are satisfied for every $n\in \N$. \begin{enumerate}
    \item[(i)] $k(\mathcal{T}_n)>0$.
    \item[(ii)] If $(U_m)$ is an increasing sequence of open sets so that $\bigcup\{U_m ~:~ m\in\N\}\in \mathcal{T}_n$ then there exists $m_0\in \N$ so that $U_{m_0}\in\mathcal{T}_n$.
    \item[(iii)] If $U\in\mathcal{T}_n$ and $V$ is an open set so that $U\Delta V$ is meagre then $V\in \mathcal{T}_n$.
\end{enumerate}
\end{definition}

Ryll-Nardzewski obtained a necessary and sufficient condition of the existence of a $\sigma$-additive, strictly positive measure on a complete Boolean algebra, see \cite[Adendum]{Kelley1959}. It is easy to see that property (***) is equivalent to $\roa{K}$ satisfying Ryll-Nardzewski's condition.  The following theorem of Ayerbe Toledano \cite{Ayerbe_Toledano1990} therefore follows immediately from Ryll-Nardzewski's result.

\begin{theorem}\label{Theorem:  Category measure characterization}
A compact Hausdorff space $K$ admits a category measure if and only if it satisfies property (***).
\end{theorem}

\begin{remark}
Cambern and Greim \cite{CambernGreim1989}, see also \cite{DalesDashiellLauStrass2016}, define a category measure to be a (not necessarily finite) positive Borel measure $\mu$ on a Stonean space $K$ which satisfies the following conditions: \begin{enumerate}
    \item[(i)] $\mu$ is regular on closed sets of finite measure,
    \item[(ii)] $\mu(N)=0$ for every meagre Borel set $N$, and,
    \item[(iii)] for every nonempty clopen set $A$ there exists a clopen set $A_0\subseteq A$ so that $0<\mu(A_0)<\infty$.
\end{enumerate}
Such measures are called \emph{perfect} in \cite{BehrendsDanckwertsEvansGobelGreimMeyfarthWinfried1977}, and are distinct from the category measures discussed here.
\end{remark}

\section{Characterizations of $\contKd$}\label{Section:  Characterisations of C(K)^d}

\subsection{The Gleason cover and the Maeda-Ogasawara theorem}\label{Subsection:  The Gleason cover and the Maeda-Ogasawara theorem}

Let $K$ be a compact Hausdorff space. A seminal result of Gleason \cite{Gleason1958} associates, in a canonical way, with $K$ a Stonean space $G_K$, its projective cover.  In order to formulate this result we recall the following.  If $L$ is compact Hausdorff space, a continuous surjection $f:K\to L$ is called \emph{irreducible} if $f[C]\neq L$ for every proper closed subset $C$ of $K$.

\begin{theorem}\label{Theorem:  Gleason Cover}
Let $K$ be a compact Hausdorff space.  There exists a Stonean space $G_K$, unique up to homeomorphism, and an irreducible map $\pi_K:G_K\to K$.
\end{theorem}

Because $\pi_K$ is onto, the induced map $T_{\pi_K}:\contK\ni f\mapsto f\circ \pi_K\in \cont (G_K)$ is an isometric linear lattice isomorphism onto a closed vector lattice subspace of $\cont(G_K)$.  Moreover, since $\pi_K$ is irreducible, $\pi_K^\ast[\contK]$ is a order dense sublattice of $\cont (G_K)$.  Clearly, $\pi_K^\ast[\contK]$ is majorising in $\cont(G_K)$, seeing as it contains the order unit $\onef_{G_K}$ of $\cont(G_K)$. Lastly we note that, $G_K$ being Stonean, $\cont(G_K)$ is Dedekind complete.  Hence we have the following.

\begin{theorem}\label{Theorem:  Dedekind Completion ito Gleason Cover}
Let $K$ be a compact Hausdorff space.  Then $\contKd$ is isometrically lattice isomorphic to $\cont(G_K)$.
\end{theorem}

Let $\bar\R$ be the two-point compactification of $\R$.  Denote by $\cont^\infty(K)$ the space of continuous functions $f:K\to\bar \R$ so that $f^{-1}[\R]$ is dense (hence open and dense) in $K$.  If $K$ is a Stonean space then $\cont^\infty(K)$ is a universally complete vector lattice, see for instance \cite[Theorem 7.27]{AliprantisBurkinshaw78}, and $\contK$ is the ideal in $\cont^\infty(K)$ generated by $\onef$.

A classical result in the representation theory for vector lattices is due to Maeda and Ogasawara \cite{MaedaOgasawara1942}, see for instance \cite{AliprantisBurkinshaw78} for a more recent presentation.

\begin{theorem}\label{Theorem:  Maeda-Ogasawara Representation}
Let $E$ be an Archimedean vector lattice with weak order unit $e$.  There exists a Stonean space $\Omega_E$ and a linear lattice isomorphism $T:E\to \cont^\infty(\Omega_E)$ onto an order dense sublattice of $\cont^\infty(\Omega_E)$ so that $Te=\onef$.
\end{theorem}

In Theorem \ref{Theorem:  Maeda-Ogasawara Representation}, let $E=\contK$.  Then $T(\onef_K)=\onef_{\Omega_E}$.  For $f\in \cont^\infty(\Omega_E)$, $f\in\cont(\Omega_E)$ if and only if there exists $c>0$ so that $|f|\leq c\onef_{\Omega_E}$.  Therefore $T[\contK]$ is an order dense and majorizing vector lattice subspace of $\cont^\infty(\Omega_E)$.  Hence we have the following.

\begin{theorem}\label{Theorem:  Dedekind Completion ito Gleason Cover}
Let $K$ be a compact Hausdorff space.  Then $\contKd$ is isometrically lattice isomorphic to $\cont(\Omega_{\contK})$.
\end{theorem}

A compact Hausdorff space is uniquely determined, up to homeomorphism, by its lattice of real-valued continuous functions. The Dedekind completion of a vector lattice is unique up to a linear lattice isomorphism. Therefore $\Omega_{\contK}$ and $G_K$ are homeomorphic.  This can be seen directly by recalling the constructions of $\Omega_{\contK}$ and $G_K$, respectively.  $G_K$ may be constructed as the Stone space of $\roa{K}$.  On the other hand, $\Omega_{\contK}$ is the Stone space of the Boolean algebra of bands in $\contK$ which is isomorphic to $\roa{K}$.

\subsection{Category measurable functions}\label{Subsection:  Category measurable functions}

Denote by $\catb{K}$ the Archimedean vector lattice consisting of all real-valued, bounded and $\ca{K}$-measurable functions on $K$.  The subset
\[
\catn{K}=\{f\in \catb{K}~:~K\setminus f^{-1}[\{0\}]\in\meagre{K}\}
\]
of $\catb{K}$ is a $\sigma$-ideal in $\catb{K}$.  Therefore
\[
\catl{K} = \catb{K}/\catn{K}
\]
is an Archimedean vector lattice \cite[Theorem 60.3]{LuxemburgZaanen1971RSI}.  For each $f\in \catb{K}$ we denote by $\hat f$ the equivalence class in $\catl{K}$ generated by $f$.  We note that for $\hat f,\hat g\in \catl{K}$,
\begin{eqnarray}
\hat f\leq \hat g \text{ if and only if } \{x\in K~:~f(x)>g(x)\}\in\meagre{K}\label{Proposition:  Characterisation order in CatM}.
\end{eqnarray}

For $\hat f\in\catl{K}$ set
\[
\|\hat f\|_\infty = \inf\{c\geq 0 ~:~ |\hat f|\leq c\hat \onef\}.
\]
Note that for $\hat f \in \catl{K}$,
\[
\|\hat f\|_\infty = \inf\{c\geq 0 ~:~ \{x\in K ~:~ |f(x)|>c\}\in\meagre{K}\}.
\]
With respect to this norm the vector lattice $\catl{K}$ is a Banach lattice.  Because $\ca{K}/\meagre{K}$ is complete, $\catl{K}$ is Dedekind complete, see for instance \cite[363M \& 363N]{Fremlin2002}.  We remark that it is possible to prove this result directly.  Every equivalence class in $\catl{K}$ contains a bounded lower semi-continuous function, and, every bounded lower semi-continuous function belongs to $\catb{K}$.  Dedekind completeness of $\catl{K}$ then follows from the fact that the pointwise supremum of any family of lower semi-continuous functions is again lower semi-continuous.

\begin{theorem}\label{Theorem:  Dedekind completion ito Linfty}
Let $K$ be a compact Hausdorff space.  Then $\contKd$ is isometrically lattice isomorphic to $\catl{K}$.
\end{theorem}
\begin{proof}

$\contK$ is vector lattice subspace of $\catb{K}$.  Furthermore, for $\contK\cap \catn{K}=\{\zerof\}$.  Therefore the mapping $T:\contK\ni u\mapsto \hat u\in\catl{K}$ is a linear lattice isomorphism onto a vector lattice subspace of $\catl{K}$.  Clearly $T$ is an isometry and $T[\contK]$ is majorising in $\catl{K}$.

Since $\catl{K}$ is a Dedekind complete Banach lattice, it remains only to verify that $T[\contK]$ is order dense in $\catl{K}$.  To this end, consider $\hat\zerof<\hat u\in\catl{K}$.  There exists a real number $\epsilon>0$, an open set $U$ in $K$ and $N\in\meagre{K}$ so that $\epsilon<u(x)$ for every $x\in U\Delta N$.  According to (\ref{Proposition:  Characterisation order in CatM}) there exists $M\in\meagre{K}$ so that $0\leq u(x)$ for all $x\in K\setminus M$. Consider a function $v\in\contK$ so that $\zerof < v\leq \epsilon \onef_U$. Fix $x\in K\setminus (M\cup N)$.  If $x\in U$ then $v(x)\leq \epsilon< u(x)$.  If $x\in K\setminus U$ then $v(x)=0\leq u(x)$.  Therefore $v(x)\leq u(x)$ for all $x\in K\setminus (M\cup N)$ so that $Tv=\hat v\leq \hat u$ by (\ref{Proposition:  Characterisation order in CatM}).

\end{proof}

\section{The dual of $\contKd$}\label{Section:  The dual of C(K)^d}

From the results discussed in Section \ref{Section:  Characterisations of C(K)^d}, in particular Theorems \ref{Theorem:  Dedekind Completion ito Gleason Cover} and \ref{Theorem:  Dedekind completion ito Linfty}, we obtain immediately two characterizations of the dual of $\contKd$.  The first follows directly from the Riesz Representation Theorem for compact Hausdorff spaces.

\begin{theorem}\label{Theorem:  Dedekind completion dual ito Gleason Cover}
Let $K$ be a compact Hausdorff space.  Then $(\contKd)^\ast$ is isometrically lattice isomorphic to $\borm(G_K)$.
\end{theorem}

Theorem \ref{Theorem:  Dedekind completion ito Linfty} yields a second characterization of $(\contKd)^\ast$.  Denote by $\cfam{K}$ the set of bounded, finitely additive signed measures $\mu$ on $\ca{K}$ with the property that $\mu(N)=0$ for every $N\in\meagre{K}$.  The space $\cfam{K}$ is a vector lattice \cite[Section 1]{YosidaHewitt}.  With respect to the total variation norm
\[
\|\varphi\| = |\varphi|(K),~ \varphi\in \cfam{K}
\]
$\cfam{K}$ is a Banach lattice, see for instance \cite{Luxemburg1991}.  We now have the following.

\begin{theorem}\label{Theorem:  Dedekind completion dual ito Category Measures}
Let $K$ be a compact Hausdorff space.  Then $(\contKd)^\ast$ is isometrically lattice isomorphic to $\cfam{K}$.  In particular, for every $\varphi\in (\contKd)^\ast$ there exists a unique $\mu_\varphi\in \cfam{K}$ so that
\[
\varphi(\hat u)=\int_K ud\mu_\varphi,~ \hat u\in\contKd,
\]
and the mapping $S:(\contKd)^\ast\ni \varphi\mapsto \mu_\varphi\in \cfam{K}$ is a linear lattice isometry onto $\cfam{K}$.
\end{theorem}

\begin{proof}
It follows immediately from \cite[Theorem 2.3]{YosidaHewitt} that $S$ is a linear isometry onto $\cfam{K}$.  That both $S$ and its inverse are positive operators follows immediately from the construction of $\mu_\varphi$ from $\varphi\in (\contKd)^\ast$ and the definition of the integral.  Indeed, for $\varphi\in (\contKd)^\ast_+$ and $\mu\in \cfam{K}_+$
\[
\mu_\varphi(B) = \varphi(\hat \onef_B),~ B\in\ca{K}
\]
and
\[
\int_K fd\mu = \sup_{s\in S_f} \int_K sd\mu,~ f\in \catb{K}_+
\]
where $S_f$ consists of all simple, positive functions dominated by $f$, see \cite{Luxemburg1991} and \cite{YosidaHewitt} for the details.
\end{proof}

Recall that $\ca{K}/\meagre{K}$ is a complete Boolean algebra, isomorphic to $\roa{K}$.  For $B\in\ca{K}$ let $\hat B$ denote the equivalence class in $\ca{K}/\meagre{K}$ containing $B$.  For $B_0,B_1\in\ca{K}$, $\hat B_0 = \hat B_1$ if and only if $B_0\Delta B_1\in\meagre{K}$.  Denote by $\famr{K}$ the space of bounded finitely additive measures on $\roa{K}$.  This space is a Banach lattice with respect to the pointwise order and variation norm, see for instance \cite[326Y (j)]{Fremlin2002}.

Let $\varphi\in\cfam{K}$. Since $\varphi[\meagre{K}]=\{0\}$, $\varphi$ induces a (signed) finitely additive measure $\mu_\varphi$ on $\roa{K}$,
\[
\mu_\varphi (\hat B) = \varphi(B),~ B\in \ca{K}.
\]
Conversely, every finitely additive measure $\mu$ on $\roa{K}$ induces a measure $\varphi^\mu\in\cfam{K}$,
\[
\varphi^\mu (B) = \mu(\hat B),~ B\in\ca{K}.
\]
The maps $\cfam{K}\ni\varphi\mapsto \mu_\varphi\in\famr{K}$ and $\famr{K}\ni \mu\mapsto \varphi^\mu\in\cfam{K}$ are positive, linear isometries, and each is the inverse of the other.  Therefore we have the following.

\begin{corollary}\label{Corollary:  Dual of C(K)^d ito measures on ROA(K)}
Let $K$ be a compact Hausdorff space.  The following statements are true. \begin{enumerate}
    \item[(i)] The following spaces are pairwise isometrically lattice isomorphic: $(\contKd)^\ast$, $\famr{K}$, $\cfam{K}$ and $\borm(G_K)$.
    \item[(iv)] If $K$ is Stonean then $\contK^\ast$, $\famr{K}$ and $\bormK$ pairwise isometrically lattice isomorphic.
\end{enumerate}
\end{corollary}

 All this is well (but perhaps not widely) known, see for instance \cite[Chapter 32]{Fremlin2002}.  That $\cfam{K}$ and $\borm(G_K)$ are isometrically lattice isomorphic can also be derived from \cite[Paragraph 4.5]{YosidaHewitt}, and \cite[362A]{Fremlin2002} may be used to show that $\famr{K}$ is isometrically lattice isomorphic to $\borm(G_K)$.

For the remainder of the paper we will be primarily concerned with the representation of $(\contKd)^\ast$ given in Theorem \ref{Theorem:  Dedekind completion dual ito Category Measures}.  

\section{Order continuous elements of $(\contKd)^\ast$ and hyper-Stonean spaces}\label{Section:  Order continuous elements of (C(K)^d)^*}

The first result of this section is a characterisation of $\sigma$-order continuous functionals on $\contKd$.  In contrast with the classical representation theorem for functionals on $\contK$ in terms of Borel measures, there is an exact correspondence between countably additivity of a measure and $\sigma$-order continuity of the corresponding functional.

\begin{theorem}\label{Theorem:  Sigma-order continuous functionals}
A measure $\varphi\in(\contKd)^\ast$ is $\sigma$-order continuous if and only if it is countably additive on $\ca{K}$.
\end{theorem}

\begin{proof}
Both the countably additive elements of $\cfam{K}$ and the $\sigma$-order continuous elements in $(\contKd)^\ast$ are ideals in the respective ambient spaces.  Therefore we may assume that $\varphi\geq 0$.

Assume that $\varphi$ is $\sigma$-order continuous.  Let $(C_n)$ be a decreasing (with respect to inclusion) sequence of sets in $\Sigma_C$ so that $\bigcap\{C_n ~:~ n\in\N\}=\emptyset$.  For each $n\in \N$ let $f_n = \onef_{C_n}$.  Then $\hat f_n\downarrow \hat \zerof$ in $\catl{K}$.  Therefore
\[
\varphi(C_n) = \int_K \hat f_n d\varphi\longrightarrow 0
\]
so that $\varphi$ is countably additive.

Conversely, assume that $\varphi$ is countably additive.  Consider a sequence $\hat f_n \downarrow \hat \zerof$ in $\catl{K}$.  Replacing each $f_n$ with $(f_1\wedge\ldots\wedge f_n)\vee \zerof$ if necessary and using (\ref{Proposition:  Characterisation order in CatM}) we may assume that the sequence $(f_n)$ in $\catb{K}$ is pointwise decreasing and bounded below by $0$. Let $f:K\ni x\mapsto \inf_{n\in\N} f_n(x)\in\R$. Then $f$ is $\ca{K}$-measurable and bounded on $K$; that is, $f\in\catb{K}$.  But $\hat \zerof \leq \hat f\leq \hat f_n$ for all $n\in N$.  Therefore $\hat f=\hat \zerof$ so that $f^{-1}[\R_+]\in\mathcal{N}$. By the Lebesgue Dominated Convergence Theorem,
\[
\varphi(\hat f_n)=\int_{K} f_n d\varphi\longrightarrow \int_K fd\varphi=0.
\]
so that $\varphi$ is $\sigma$-order continuous.
\end{proof}

For a measure $\varphi\in \cfam{K}_+$ we define regularity in the same way as for Borel measures.  That is, $\varphi$ is \emph{regular} if for every $B\in\bor{K}$,
\[
\sup\{\mu(C) ~:~ C\subseteq B \text{ is compact}\} = \mu(B) = \inf\{\mu(U) ~:~ U\supseteq B \text{ is open}\}.
\]
We note that each of the identities above implies the other.  A general measure $\varphi\in \cfam{K}$ is regular if $|\varphi|$ is regular.

\begin{theorem}\label{Theorem:  Order continuous functionals}
A measure $\varphi\in(\contKd)^\ast$ is order continuous if and only if it is countably additive and regular on $\ca{K}$.
\end{theorem}

\begin{proof}
Let $\varphi\geq 0$.  Assume that $\varphi$ is order continuous. By Theorem \ref{Theorem:  Sigma-order continuous functionals}, $\varphi$ is countably additive.  Let $A\in\ca{K}$.  Then $\hat \onef_A\in\catl{K}$. Fix $\epsilon>0$ and let $D_\epsilon=\{\hat f\in \catl{K} ~:~ f\in\contK,~ (1+\epsilon/2)\onef_A\leq f\}$. According to Theorem \ref{Theorem:  Dedekind completion ito Linfty}, $D_\epsilon\downarrow (1+\epsilon/2)\hat \onef_A$ in $\catl{K}$.  By order continuity of $\varphi$,
\[
(1+\epsilon/2)\varphi(A)=\int (1+\epsilon/2)\onef_A d\phi = \inf_{\hat f\in D_\epsilon}\int_K f d\varphi.
\]
Therefore there exists $\hat f\in D_\epsilon$ so that
\[
\int_K fd\varphi < (1+\epsilon)\varphi(A).
\]
Let $U=f^{-1}[(1,\infty)]$.  Then $U\supseteq A$ is open and
\[
\varphi(U)=\int_K \onef_U d\varphi \leq \int_K fd\varphi < (1+\epsilon)\varphi(A).
\]
Therefore $\varphi(A)\leq \inf \{\varphi(U) ~:~ U\supseteq A ~\text{open}\} <(1+\epsilon)\varphi(A)$.  This holds for all $\epsilon>0$ so that $\varphi(A)= \inf \{\varphi(U) ~:~ U\supseteq A ~\text{open}\}$.  Therefore $\varphi$ is outer regular, hence regular.

Conversely, assume that $\varphi$ is countably additive and regular on $\ca{K}$.  Then the restriction $\varphi_0$ of $\varphi$ to $\bor{K}$ is countably additive and regular, and $\varphi_0(N)=0$ for every meagre Borel set $N$. Therefore $\varphi_0$ is a normal Borel measure on $K$ so that the restriction of the functional $\varphi$ to $\contK$ is order continuous.  By Lemma \ref{Lemma:  Restriction and extension of order continuous functionals} (i), $\varphi$ is order continuous on $\contKd$.
\end{proof}

\begin{proposition}\label{Proposition:  Strictly positive functionals}
Let $\varphi \in (\contKd)^\ast_+$.  Then $\varphi$ is strictly positive if and only if, for every  every $A\in\ca{K}$, $\varphi(A)=0$ if and only if $A\in\meagre{K}$.
\end{proposition}

\begin{proof}
Assume that $\varphi$ is strictly positive.  Let $A\in\ca{K}\setminus \meagre{K}$.  Then $\hat \zerof<\hat\onef_A$ so that $\varphi(A)=\varphi(\hat \onef_A)>0$.

Assume that for every $A\in\ca{K}$, $\varphi(A)=0$ if and only if $A\in\meagre{K}$.  Consider any $\hat\zerof <\hat f\in\catl{K}$.  Then there exists an $\epsilon>0$ so that $A=f^{-1}[[\epsilon,\infty]]\in\ca{K}\setminus\meagre{K}$.  Consequently,
\[
\varphi(\hat f)\geq \int_K\onef_A d\varphi = \varphi(A)>0.
\]
Therefore $\varphi$ is strictly positive.
\end{proof}

Theorems \ref{Theorem:  Category measure characterization}, \ref{Theorem:  Sigma-order continuous functionals} and \ref{Theorem:  Order continuous functionals}, and Proposition \ref{Proposition:  Strictly positive functionals} now yields the following result, establishing the relationship between order continuous functionals and category measures.

\begin{corollary}\label{Corollary:  Characterizations of strictly positive order continuous functionals}
Let $K$ be a compact Hausdorff space.  Then the following statements are equivalent. \begin{enumerate}
    \item[(i)] $K$ satisfies property (***).
    \item[(ii)] $K$ admits a category measure.
    \item[(iii)] $\contKd$ admits a strictly positive $\sigma$-order continuous linear functional.
    \item[(iv)] $\contK$ admits a strictly positive $\sigma$-order continuous linear functional.
    \item[(v)] $\contK$ admits a strictly positive order continuous linear functional.
    \item[(vi)] $\contKd$ admits a strictly positive order continuous linear functional.
\end{enumerate}
\end{corollary}

\begin{proof}
That (i) and (ii) are equivalent is Theorem \ref{Theorem:  Category measure characterization}. The equivalence of (ii) and (iii) follows immediately from the definition of a category measure, Theorem \ref{Theorem:  Sigma-order continuous functionals} and Proposition \ref{Proposition:  Strictly positive functionals}.  That (iii) implies (iv) follows from Lemma \ref{Lemma:  Restriction and extension of order continuous functionals} (ii).

To see that (iv) implies (v), assume that $\contK$ admits a strictly positive $\sigma$-order continuous linear functional $\varphi$.  Then there exists a fully supported regular Borel measure on $K$, hence $K$ satisfies cc, see for instance \cite[Proposition 4.1.6]{DalesDashiellLauStrass2016}.  Then $\contK$ is order separable so that $\varphi$ is order continuous.

The equivalence of (v) and (vi) follows from Lemma \ref{Lemma:  Restriction and extension of order continuous functionals} (i) and (ii).  Lastly, that (vi) implies (iii) is obvious, which complete the proof.
\end{proof}

\begin{remark}
For the equivalence of (iii) to (vi) in Corollary \ref{Corollary:  Characterizations of strictly positive order continuous functionals} the assumption of strict positivity is essential.  Indeed, there exists a compact Hausdorff space $K$ and a $\sigma$-order continuous functional on $\contK$ which is not order continuous \cite[Example 4.7.16]{DalesDashiellLauStrass2016}.  Furthermore, the statement `there exists a Stonean space $K$ and a $\sigma$-order continuous functional on $\contK$ which is not order continuous' is equivalent to the existence of a measurable cardinal, see \cite{Luxemburg1967}.
\end{remark}


As an application of Corollary \ref{Corollary:  Characterizations of strictly positive order continuous functionals} we obtain the following topological characterization of hyper-Stonean spaces.

\begin{theorem}\label{Theorem:  Characterisation of hyper-Stonean spaces}
Let $K$ be a Stonean space.  Then $K$ is hyper-Stonean if and only if there exists a collection $\mathcal{U}$ of clopen subsets of $K$ so that $\bigcup \mathcal{U}$ is dense in $K$ and every $U\in\mathcal{U}$ satisfies property (***).
\end{theorem}

\begin{proof}
Assume that $K$ is hyper-Stonean.  For each normal measure $\mu$ on $K$, let $S_\mu$ denote the support of $\mu$.  Let $\mathcal{U}=\{S_\mu ~:~ \mu \in\mathrm{N}(K)\}$.  By definition of a hyper-Stonean space $\bigcup\mathcal{U}$ is dense in $K$.  Each $S_\mu\in \mathcal{U}$ is clopen, hence itself Stonean, and $\mu$ defines a strictly positive order continuous functional on $\cont (S_\mu)$.  Corollary \ref{Corollary:  Characterizations of strictly positive order continuous functionals} implies that each $S_\mu$ satisfies property (***).

Assume that there exists a collection $\mathcal{U}$ of clopen subsets of $K$ so that $\bigcup \mathcal{U}$ is dense in $K$ and every $U\in\mathcal{U}$ satisfies property (***).  We claim that each $U\in \mathcal{U}$ is the support of a normal measure on $K$. Fix $U\in \mathcal{U}$.  By Corollary \ref{Corollary:  Characterizations of strictly positive order continuous functionals} there exists a strictly positive order continuous functional $\varphi$ on $\cont (U)$.  Consider the linear functional $\psi:\contK \ni f\mapsto \psi(f|U)\in \R$.  Because $U$ is clopen, hence regular closed in $K$, $\psi$ is order continuous by \cite[Theorem 3.4]{KandicVavpeticPositivity2019}.  Therefore there exists a unique normal measure $\mu$ on $K$ so that
\[
\psi(u) = \int_K fd\mu,~ f\in\contK.
\]
For every $f\in\contK_+$, $f|S_\mu=\zerof$ if and only if $\psi(f)=0$, if and only if $f|U=\zerof$.  Therefore $U=S_\mu$, which verifies our claim.

Since $\bigcup \mathcal{U}$ is dense in $K$, the union of the supports of the normal measures on $K$ is dense in $K$; that is, $K$ is hyper-Stonean.
\end{proof}

\begin{remark}
It is possible to obtain our main results, namely Theorem \ref{Theorem:  Dedekind completion dual ito Category Measures}, from which the results in Section \ref{Section:  Order continuous elements of (C(K)^d)^*} follow in a less direct manner.  We briefly recall how this may be achieved.

In \cite{deJongevanRooijRieszSpaces1977}, de Jonge and van Rooij give a construction of $\contKd$ in terms of Borel measurable functions which is very similar to that given in Theorem \ref{Theorem:  Dedekind completion ito Linfty}.  If $\mathrm B$ denotes the vector lattice of bounded Borel measurable functions on $K$ and $\mathrm N$ the subspace of $\mathrm B$ consisting of these functions which vanish of a meagre Borel set, then $\contKd$ can be identified with $\mathrm D (K)= \mathrm B/\mathrm N$.  Dales et al. \cite{DalesDashiellLauStrass2016} call the space $\mathrm D (K)$ the Dixmier algebra of $K$.  Using this construction our results can be obtained via the machinery of measures on Boolean algebras as set out, for instance, in \cite{Fremlin2002}.  We have opted for a direct and more transparent approach.
\end{remark}





\bibliographystyle{amsplain}
\bibliography{Dedekindcompletiondual}

\end{document}